\documentclass{elsart}

\usepackage{amssymb,amsfonts,amsxtra}
\usepackage{epic,eepic}
\usepackage[all]{xypic}
\xyoption{dvips}

\newcommand{\lead}{\mathrm{lead}}
\newcommand{\lexp}{\mathrm{lexp}}
\newcommand{\supp}{\mathrm{supp}}
\newcommand{\NF}{\mathrm{NF}}

\begin{document}

\begin{frontmatter}

\title{A normal form algorithm\\for the Brieskorn lattice}
\author{Mathias Schulze}
\address{Department of Mathematics\\University of Kaiserslautern\\67663 Kaiserslautern, Germany}
\ead{mschulze@mathematik.uni-kl.de}

\begin{abstract}
This article describes a normal form algorithm for the Brieskorn lattice of an isolated hypersurface singularity.
It is the basis of efficient algorithms to compute the Bernstein--Sato polynomial, the complex monodromy, and Hodge--theoretic invariants of the singularity such as the spectral pairs and good bases of the Brieskorn lattice.
The algorithm is a variant of Buchberger's normal form algorithm for power series rings using the idea of partial standard bases and adic convergence replacing termination.
\end{abstract}

\begin{keyword}
hypersurface singularity \sep Brieskorn lattice \sep Bernstein-Sato polynomial \sep monodromy \sep spectral pairs \sep good basis \sep mixed Hodge structure \sep standard basis
\MSC 13N10 \sep 13P10 \sep 32S35 \sep 32S40
\end{keyword}

\end{frontmatter}

\section{Introduction}

Isolated hypersurface singularities form the simplest class of singularities.
Their intensive study in the past has led to a variety of invariants.
The Milnor number is one of the simplest, and can easily be computed using standard basis methods.
A finer invariant is the monodromy of the singularity.
E.~Brieskorn \cite{Bri70} developed the theoretical background for computing the complex monodromy.
He gave an ad hoc definition of an object $H''$, later called the Brieskorn lattice.
Its great importance was a priori not clear.
The complex monodromy can be expressed in terms of the differential structure of the Brieskorn lattice.

The finest known invariants come from a mixed Hodge structure associated to an isolated hypersurface singularity.
The notion of a mixed Hodge structure was introduced by P.~Deligne \cite{Del70a} as a generalization of the classical Hodge structure on the cohomology of a compact K\"ahler manifold. 
J.H.M.~Steenbrink \cite{Ste76} defined this mixed Hodge structure in terms of resolutions of singularities.
A.N.~Varchenko \cite{Var82a} and later J.~Scherk and J.H.M.~Steenbrink \cite{SS85} described this mixed Hodge structure in terms of the differential structure of the Brieskorn lattice. 
The mixed Hodge numbers correspond to the spectral pairs and determine the complex monodromy.
The spectral pairs have a semicontinuity property \cite{Ste85} with respect to unfoldings of the singularity.

Based on properties of the mixed Hodge structure, M.~Saito \cite{Sai89} constructed two endomorphisms $A_0$ and $A_1$ of the Milnor algebra.
These two endomorphisms determine the differential structure of the Brieskorn lattice and, immediately, the above invariants.

The Bernstein--Sato polynomial is associated to a general complex polynomial \cite{Ber72} or convergent power series \cite{Bjo72}.
T.~Oaku \cite{Oak97} presented the first algorithm to compute it in the global case.
A new method by M.~Noro \cite{Nor02} is impressively faster.
In the isolated singularity case, B.~Malgrange \cite{Mal75} described the Bernstein--Sato polynomial in terms of the Brieskorn lattice and its close relation to the complex monodromy.

In \cite{Sch01,SS01,Sch02a,Sch02c}, we have developed algorithmic methods to compute all of the above invariants of isolated hypersurface singularities.
There is an implementation \cite{Sch02b,Sch03c} of these algorithms in the computer algebra system {\sc Singular} \cite{GPS03}.
Our algorithm to compute the complex monodromy is much faster and can compute much more difficult examples than Brieskorn's algorithm.
Our algorithm to compute the local Bernstein--Sato polynomial is based on B.~Malgrange's description in terms of the Brieskorn lattice.
It is much faster than M.~Noro's algorithm since computations in rings of differential operators are replaced by computing their action on power series rings.
However, it is restricted to the isolated singularity case while M.~Noro's algorithm works in the general global case.
All our algorithms require the computation of a basis representation in the Brieskorn lattice.
In \cite[Sec.~10.2]{Sch01} this is done by a sequence of divisions by the Jacobian ideal which is, in general, very hard to compute.
The subject of this article is a normal form algorithm for the Brieskorn lattice replacing this sequence of full divisions by a sequence of partial divisions.
This new method turns out to be much more efficient.

In the first section, we recall the definition and the main properties of the Brieskorn lattice.
We introduce the formal Brieskorn lattice and describe it as a cokernel of a formal family of differential operators which is finite over the base.

In the second section, we consider such a formal family of differential operators in general.
We describe a normal form algorithm to compute a presentation of the cokernel which is a finitely generated module over the formal power series ring in the parameters of the family.
This algorithm is a variant of B.~Buchberger's \cite{Buc65,Buc85} normal form algorithm.
There are three major differences compared to the classical algorithm:
\begin{enumerate}
\item The polynomial ring is replaced by a formal power series ring.
Termination of the algorithm is replaced by adic convergence.
\item The standard basis is replaced by a partial standard basis, a set of power series which specializes to a standard basis.
\item There is only a module structure with respect to the parameters of the family and the partial standard basis is not finite.
\end{enumerate}
Although the algorithm does not terminate in general, it serves to compute exact results by using appropriate degree bounds.
The algorithms in \cite{Sch01,Sch02c} implicitly contain such degree bounds to compute the above invariants of isolated hypersurface singularities.
There are also a priori degree bounds in \cite{Sch02a}, but they are useless in practice.
Essentially the double number of variables plus the double Milnor number is a degree bound that satisfies all requirements.

In the third section, we demonstrate the power of our algorithm.
We apply the {\sc Singular} implementation \cite{Sch03c} to examples from \cite{Sch01,Sch02c} and \cite{Nor02} and list the timings.

Families which are finite over the base occur in many situations in algebraic geometry and singularity theory.
For example, A.~Fr\"uhbis--Kr\"uger \cite{FK02} has developed algorithms to compute moduli spaces and adjacencies of singularities based on the idea of partial standard bases.
One can expect more applications of our methods in the future.

We shall denote row vectors by a lower bar, column vectors by an upper bar, row indices by lower indices, and column indices by upper indices.

\begin{ack}
I express my gratitude to A.~Fr\"uhbis--Kr\"uger and G.--M.~Greuel for fruitful discussions and to the two referees for useful hints to improve this article.
I wish to thank M.~Granger for pointing out an error in the preprint version and P.~Bitsch for proofreading my English.
\end{ack}

\section{The formal Brieskorn lattice}

Let $\SelectTips{cm}{}\xymatrix{f:U\ar[r]&\Cset}$ be a holomorphic function on an open neighbourhood $0\in U\subset\Cset^n$ of the origin.
We choose a system of complex coordinates $\underline x=x_1,\dots,x_n$ at $0\in\Cset^n$ and denote by $\underline\partial=\partial_1,\dots,\partial_n=\partial_{x_1},\dots,\partial_{x_n}$ the corresponding derivatives such that the commutator of $\partial_i$ and $x_j$ is $[\partial_i,x_j]=\delta_{i,j}$.
We consider $f$ as a germ of a holomorphic function at $0\in\Cset^n$, which means that $U$ can be arbitrarily small.
This is equivalent to considering the convergent power series $f\in\Cset\{\underline x\}$.
We assume that $f(0)=0$, and that the origin is an isolated critical point of $f$.
This means that $0\in U$ is the only point with $\partial_1(f)(0)=\cdots=\partial_n(f)(0)=0$ for some $U$, or, more algebraically, that $\langle\underline x\rangle^m\subset\langle\underline\partial(f)\rangle\subset\langle\underline x\rangle$ for some $m\ge1$.
The complex dimension 
\[
\mu=\dim_\Cset\bigl(\Cset\{\underline x\}/\langle\underline\partial(f)\rangle\bigr)<\infty
\]
is called the Milnor number.
By the finite determinacy theorem \cite[Thm.~9.1.4]{dJP00}, one can choose, in this case, a coordinate system $\underline x$ such that $f\in\Cset[\underline x]$ is a polynomial.
 
We denote by $\Omega^\bullet=\Omega^\bullet_{\Cset^n,0}$ the complex of germs of holomorphic differential forms at $0\in\Cset^n$.
Its elements are differential forms with coefficients in the convergent power series ring $\Cset\{\underline x\}$.
The Brieskorn lattice \cite{Bri70} is defined by
\[
H''=\Omega^n/\d f\wedge\d\Omega^{n-2}
\]
and becomes a $\Cset\{t\}$--module by setting
\begin{equation}\label{22}
t\cdot[\omega]=[f\omega]
\end{equation}
for $[\omega]\in H''$.
By M.~Sebastiani \cite{Seb70}, $H''$ is a free $\Cset\{t\}$--module of rank $\mu$.
We denote by $\Omega$ the $\mu$--dimensional $\Cset$--vector space
\[
\Omega=\Omega^n/\d f\wedge\Omega^{n-1}\cong\Cset\{\underline x\}/\langle\underline\partial(f)\rangle.
\]
The operators $\d$ and $\d f=\d f\wedge\cdot$ define two exact sequences.

\begin{lem}[Poincar\'e Lemma]\label{1}
\[
\SelectTips{cm}{}\xymatrix{
0\ar[r] & \Cset\ar[r] & \Cset\{\underline x\}\ar[r]^-\d & \Omega^1\ar[r]^-\d & \cdots\ar[r]^-\d & \Omega^n\ar[r] & 0
}
\]
is an exact sequence of $\Cset$--vector spaces.
\end{lem}

Since completion is exact, Lemma \ref{1} remains valid when replacing $\Omega^\bullet$ by its $\langle\underline x\rangle$--adic completion $\widehat\Omega^\bullet$.
The elements of $\widehat\Omega^\bullet$ are differential forms with coefficients in the formal power series ring $\Cset[\![\underline x]\!]$.

\begin{lem}[De Rham Lemma]\label{2}
\[
\SelectTips{cm}{}\xymatrix{
0\ar[r] & \Cset\{\underline x\}\ar[r]^-{\d f} & \Omega^1\ar[r]^-{\d f} & \cdots\ar[r]^-{\d f} & \Omega^n\ar[r] & \Omega\ar[r] & 0
}
\]
is an exact sequence of $\Cset\{\underline x\}$--modules.
\end{lem}

Also Lemma \ref{2} remains valid when replacing $\Omega^\bullet$ by $\widehat\Omega^\bullet$.
From Lemma \ref{1} and \ref{2} follows that one can define a $\Cset$--linear operator $s$ on $H''$ by
\begin{equation}\label{23}
s\cdot[\d\eta]=[\d f\wedge\eta]
\end{equation}
for $[\d\eta]\in H''$.
From Lemma \ref{1} follows that $s$ is injective.
The image of $s$ is $sH''=\d f\wedge\Omega^n/\d f\wedge\d\Omega^{n-2}$ and hence
\[
H''/sH''=\Omega.
\]
Also $s$ defines a module structure on $H''$ over a power series ring.
This power series ring is the ring
\[
\Cset\{\!\{s\}\!\}=\Bigl\{\sum_{i=0}^\infty a_is^i\in\Cset[\![s]\!]\;\Big\vert\sum_{i=0}^\infty\frac{a_i}{i!}t^i\in\Cset\{t\}\Bigr\}\subset\Cset[\![s]\!]
\]
of microlocal operators with constant coefficients and, by F.~Pham \cite{Pha77}, $H''$ is a free $\Cset\{\!\{s\}\!\}$--module of rank $\mu$.
From the definitions of $t$ and $s$ follows immediately that the commutator of $t$ and $s$ is
\[
[t,s]=s^2.
\]
We define a $\Cset$--linear operator $\partial_s$ on the localization $H''\otimes_{\Cset\{\!\{s\}\!\}}\Cset\{\!\{s\}\!\}[s^{-1}]$ by 
\begin{equation}\label{24}
t=s^2\partial_s.
\end{equation}
Then $t$ is a differential operator on $H''$ with respect to the $\Cset\{\!\{s\}\!\}$--structure.
There is also $\Cset$--linear operator $\partial_t$ on the localization $H''\otimes_{\Cset\{t\}}\Cset\{t\}[t^{-1}]$ defined by 
\begin{equation}\label{25}
s=\partial_t^{-1}.
\end{equation}
Then the commutator of $\partial_t$ and $t$ is $[\partial_t,t]=1$ and hence $\partial_t$ is a derivative by $t$.

\begin{defn}\
\begin{enumerate}
\item We call the topology induced by the $\langle\underline x\rangle$--adic topology on $\Omega^n$ on the quotient $H''$ the $\langle\underline x\rangle$--adic topology on $H''$.
\item We call the completion $\widehat H''$ of $H''$ with respect to the $\langle\underline x\rangle$--adic topology the formal Brieskorn lattice.
\end{enumerate}
\end{defn}

\begin{lem}\label{26}
There is a natural isomorphism
\[
\widehat H''=\widehat\Omega^n/\d f\wedge\d\widehat\Omega^{n-2}.
\]
\end{lem}
\begin{pf}
By definition,
\[
\widehat H''=\underset{\underset{k}{\longleftarrow}}{\lim}\bigl(\Omega^n\big/\bigl(\langle\underline x\rangle^k\Omega^n+\d f\wedge\d\Omega^{n-2}\bigr)\bigr).
\]
Since $\underline0\in\Cset^n$ is a critical point of $f$, $\langle\underline\partial(f)\rangle\subset\langle\underline x\rangle$ and hence
\begin{align*}
\d f\wedge\d\bigl(\langle\underline x\rangle^k\Omega^{n-2}\bigr)&\subset\langle\underline x\rangle^k\Omega^n,\\
\d f\wedge\d\bigl(\langle\underline x\rangle^k\widehat\Omega^{n-2}\bigr)&\subset\langle\underline x\rangle^k\widehat\Omega^n.
\end{align*}
Since $\Omega^n/\langle\underline x\rangle^k\Omega^n=\widehat\Omega^n/\langle\underline x\rangle^k\widehat\Omega^n$, this implies that
\[
\Omega^n\big/\bigl(\langle\underline x\rangle^k\Omega^n+\d f\wedge\d\Omega^{n-2}\bigr)=\widehat\Omega^n\big/\bigl(\langle\underline x\rangle^k\widehat\Omega^n+\d f\wedge\d\widehat\Omega^{n-2}\bigr)
\]
defines a natural isomorphism of inverse systems.
Hence,
\[
\widehat H''=\underset{\underset{k}{\longleftarrow}}{\lim}\bigl(\widehat\Omega^n\big/\bigl(\langle\underline x\rangle^k\Omega^n+\d f\wedge\d\widehat\Omega^{n-2}\bigr)\bigr)=\widehat\Omega^n/\d f\wedge\d\widehat\Omega^{n-2}.\qed
\]
\end{pf}

The following theorem \cite[Prop.~3.3]{Bri70} is essential for Brieskorn's algorithm to compute the complex monodromy, which is based on the $\Cset\{t\}$--structure of the Brieskorn lattice.

\begin{thm}\label{39}
The $\langle t\rangle$--adic and $\langle\underline x\rangle$--adic topology on $H''$ coincide.
In particular, the $\langle t\rangle$--adic completion of $H''$ is naturally isomorphic to $\widehat H''$ and $\widehat H''$ is a free $\Cset[\![t]\!]$--module of rank $\mu$.
\end{thm}

The $\Cset\{\!\{s\}\!\}$--structure of the Brieskorn lattice is more algebraic and, therefore, more appropriate for computational purposes.
The following proposition \cite[Prop.~7]{Sch02c} is the analogue of theorem \ref{39} for the $\Cset\{\!\{s\}\!\}$--structure, but it is much easier to prove.

\begin{prop}\label{3}
The $\langle s\rangle$--adic and $\langle\underline x\rangle$--adic topology on $H''$ coincide.
In particular, the $\langle s\rangle$--adic completion of $H''$ is naturally isomorphic to $\widehat H''$ and $\widehat H''$ is a free $\Cset[\![s]\!]$--module of rank $\mu$.
\end{prop}
\begin{pf}
We denote
\begin{align*}
\d\underline x&=\d x_1\wedge\cdots\wedge\d x_n,\\
\d\underline x_{\widehat i}&=\d x_1\wedge\cdots\wedge\d x_{i-1}\wedge\d x_{i+1}\wedge\cdots\wedge\d x_n
\end{align*}
for $1\le i\le n$.
Let
\[
\bigl[g\partial_i(f)\d\underline x\bigr]\in\bigl(\langle\underline\partial(f)\rangle^{2k}\d\underline x+\d f\wedge\d\Omega^{n-2}\bigr)\big/\d f\wedge\d\Omega^{n-2}\subset H''
\]
for some $k\ge1$.
By (\ref{23}),
\begin{align*}
\bigl[g\partial_i(f)\d\underline x\bigr]
&=\bigl[(-1)^{i+1}\d f\wedge\bigl(g\d\underline x_{\widehat i}\bigr)\bigr]\\
&=s\bigl[(-1)^{i+1}\d\bigl(g\d\underline x_{\widehat i}\bigr)\bigr]\\
&=s\bigl[\partial_i(g)\d\underline x\bigr]\\
&\in s\bigl(\bigl(\langle\underline\partial(f)\rangle^{2(k-1)}\d\underline x+\d f\wedge\d\Omega^{n-2}\bigr)\big/\d f\wedge\d\Omega^{n-2}\bigr)
\end{align*}
and hence, by induction,
\[
\bigl(\langle\underline\partial(f)\rangle^{2k}\d\underline x+\d f\wedge\d\Omega^{n-2}\bigr)\big/\d f\wedge\d\Omega^{n-2}\subset s^kH''.
\]
Since $\underline0\in\Cset^n$ is an isolated critical point of $f$, $\langle\underline x\rangle^m\subset\langle\underline\partial(f)\rangle\subset\langle\underline x\rangle$ for some $m\ge1$ and hence
\begin{multline*}
\bigl(\langle\underline x\rangle^{2km}\d\underline x+\d f\wedge\d\Omega^{n-2}\bigr)\big/\d f\wedge\d\Omega^{n-2}\\
\subset\bigl(\langle\underline\partial(f)\rangle^{2k}\d\underline x+\d f\wedge\d\Omega^{n-2}\bigr)\big/\d f\wedge\d\Omega^{n-2}.
\end{multline*}
This implies that
\begin{multline*}
\bigl(\langle\underline x\rangle^{2km}\d\underline x+\d f\wedge\d\Omega^{n-2}\bigr)\big/\d f\wedge\d\Omega^{n-2}\subset s^kH''\\
\subset\bigl(\langle\underline x\rangle^k\d\underline x+\d f\wedge\d\Omega^{n-2}\bigr)\big/\d f\wedge\d\Omega^{n-2}.
\end{multline*}
Hence, the $\langle s\rangle$--adic and $\langle\underline x\rangle$--adic topology on $H''$ coincide.\qed
\end{pf}

Note that the formal Brieskorn lattice is a $(t,s)$--module in the sense of D.~Barlet \cite{Bar93,Bar00}.
The following proposition \cite[Prop.~8]{Sch02c} describes the $\Cset[\![s]\!]$--module $\widehat H''$ as a quotient of the power series ring $\Cset[\![s,\underline x]\!]$.
It will lead to a normal form algorithm for $\widehat H''$ in the next section.

\begin{prop}\label{4}
$\d\underline x$ induces a $\Cset[\![s]\!]$--isomorphism
\[
\SelectTips{cm}{}\xymatrix{
\widehat H''=\widehat\Omega^n[\![s]\!]/(\d f-s\d)\widehat\Omega^{n-1}[\![s]\!] & \Cset[\![s,\underline x]\!]/\langle\underline\partial(f)-s\underline\partial\rangle\Cset[\![s,\underline x]\!]\ar[l]^-\sim_-{\d\underline x}.
}
\]
\end{prop}
\begin{pf}
Since
\[
\d f\wedge\d\widehat\Omega^{n-2}=(\d f-s\d)\d\widehat\Omega^{n-2}\subset(\d f-s\d)\widehat\Omega^{n-1}[\![s]\!]
\]
and by Lemma \ref{26} and (\ref{23}), there is a natural $\Cset[\![s]\!]$--linear map 
\[
\SelectTips{cm}{}\xymatrix{\widehat H''\ar[r]^-\phi & \widehat\Omega^n[\![s]\!]/(\d f-s\d)\widehat\Omega^{n-1}[\![s]\!]}.
\]
Let $\omega=\sum_{k\ge0}\omega_k s^k\in\widehat\Omega^{n-1}[\![s]\!]$ with $(\d f-s\d)\omega\in\widehat\Omega^n$.
Then $\d f\wedge\omega_{k+1}=\d\omega_k$ and hence, by (\ref{23}),
\[
s[\d\omega_{k+1}]=[\d f\wedge\omega_{k+1}]=[\d\omega_k]\in\widehat H''
\]
for all $k\ge0$.
In particular, $[\d\omega_0]\in\bigcap_{k\ge0}s^k\widehat H''=\{0\}$ and hence, by Lemma \ref{26},
\[
\d\omega_0\in\d f\wedge\d\widehat\Omega^{n-2}=\d(\d f\wedge\widehat\Omega^{n-2}).
\]
By Lemma \ref{1}, this implies that $\omega_0\in\d\widehat\Omega^{n-2}+\d f\wedge\widehat\Omega^{n-2}$ and hence
\[
(\d f-s\d)\omega=\d f\wedge\omega_0\in\d f\wedge\d\widehat\Omega^{n-2}.
\]
This shows that
\[
(\d f-s\d)\widehat\Omega^{n-1}[\![s]\!]\cap\widehat\Omega^n=\d f\wedge\d\widehat\Omega^{n-2}
\]
and hence, by Lemma \ref{26}, that $\phi$ is injective.
By Lemma \ref{1}, $\d\widehat\Omega^{n-1}=\widehat\Omega^n$, and hence $\phi$ is surjective.

For $\eta=\sum_{i=1}^n(-1)^{i+1}g_i\d\underline x_{\widehat i}\in\widehat\Omega^{n-1}[\![s]\!]$,
\[
(\d f-s\d)\eta=\sum_{i=1}^n(\partial_i(f)g_i-s\partial_i(g_i))\d\underline x=(\underline\partial(f)-s\underline\partial)\overline g\d\underline x.
\]
Hence, $\d\underline x$ induces a $\Cset[\![s]\!]$--isomorphism
\[
\SelectTips{cm}{}\xymatrix{\widehat\Omega^n[\![s]\!]/(\d f-s\d)\widehat\Omega^{n-1}[\![s]\!] & \Cset[\![s,\underline x]\!]/\langle\underline\partial(f)-s\underline\partial\rangle\Cset[\![s,\underline x]\!]\ar[l]^-\sim_-{\d\underline x}}.\qed
\]
\end{pf}

Proposition \ref{4} is the starting point for more general considerations in the next section.

\section{Formal families of differential operators}

Let $\mathbb{K}$ be a computable field and $\underline F=F_1,\dots,F_r\in\mathbb{K}[\![\underline s,\underline x]\!]\langle\underline\partial\rangle$ a formal family of differential operators where $\underline x=x_1,\dots,x_n$, $\underline\partial=\partial_1,\dots,\partial_n=\partial_{x_1},\dots,\partial_{x_n}$, and $\underline s=s_1,\dots,s_m$.
Note that the elements of $\mathbb{K}[\![\underline s,\underline x]\!]\langle\underline\partial\rangle$ are polynomial in $\underline\partial$.
The brackets $\langle\cdot\rangle$ indicate that the commutator $[x_i,\partial_i]=\delta_{i,j}$ is not zero.
We want to compute the cokernel $H=\mathbb{K}[\![\underline s,\underline x]\!]/\langle\underline F\rangle\mathbb{K}[\![\underline s,\underline x]\!]$ of the $\mathbb{K}[\![\underline s]\!]$--linear map
\[
\SelectTips{cm}{}\xymatrix{\mathbb{K}[\![\underline s,\underline x]\!]^r\ar[r]^-{\underline F}&\mathbb{K}[\![\underline s,\underline x]\!]\ar[r]\ar[r]^-{\pi_H}&H\ar[r]&0}.
\]
We assume that the specialization
\[
\underline f=f_1,\dots,f_r=\underline F(\underline s=\underline0)\in\mathbb{K}[\![\underline x]\!]
\]
is independent of $\underline\partial$ and that $\langle\underline x\rangle^k\mathbb{K}[\![\underline x]\!]\subset\langle\underline f\rangle\mathbb{K}[\![\underline x]\!]$ for some $k\ge0$.
In particular,
\[
\mu=\dim_\mathbb{K}\bigl(\mathbb{K}[\![\underline x]\!]/\langle\underline f\rangle\mathbb{K}[\![\underline x]\!]\bigr)=\dim_\mathbb{K}\bigl(H/\langle\underline s\rangle H\bigr)<\infty
\]
and hence $H$ is a finitely generated $\mathbb{K}[\![s]\!]$--module.
Then there is a matrix $D=(\underline d^j)_j\in\mathbb{K}[\![\underline s,\underline x]\!]\langle\underline\partial\rangle^{m\times r}$ such that
\[
\underline F=\underline f-\underline sD.
\]

Our considerations are motivated by the following special case.

\begin{rem}\label{9}
By Proposition \ref{4}, for $\mathbb{K}=\Cset$, $m=1$, $r=n$, $\underline f=\underline\partial(f)$, and $D=\underline\partial$,
\[
H\cong_{\Cset[\![s]\!]}\widehat H''
\]
is the formal Brieskorn lattice.
\end{rem}

Let $<_{\underline x}$ be a local degree ordering with respect to a weighted degree $\deg_{\underline x}$ on the set of monomials $\{\underline x^{\underline\alpha}\,\vert\,\underline\alpha\in\Nset^n\}=\Nset^n$ of $\mathbb{K}[\![\underline x]\!]$.
This means that
\[
\SelectTips{cm}{}\xymatrix{
\Nset^n\ar[r]^-{\deg_{\underline x}} & \Qset
}
\]
is a semigroup homomorphism with $\deg_{\underline x}(x_i)<0$, and that $<_{\underline x}$ is a semigroup ordering such that
\[
\deg_{\underline x}(\underline x^{\underline\alpha})<\deg_{\underline x}(\underline x^{\underline\beta})\Rightarrow\underline x^{\underline\alpha}<_{\underline x}\underline x^{\underline\beta}.
\]
The support of $p=\sum_{\underline\alpha} p_{\underline\alpha}\underline x^{\underline\alpha}\in\mathbb{K}[\![\underline x]\!]$ is defined by $\supp(p)=\{\underline\alpha\in\Nset^n\,\vert\,p_{\underline\alpha}\ne0\}$.
We denote the leading exponent, resp. leading term, with respect to $<_{\underline x}$ by $\lexp$, resp. $\lead$.
This means that
\begin{align*}
\lexp(p)&=\max_{<_{\underline x}}\supp(p),\\
\lead(p)&=p_{\lexp(p)}\underline x^{\lexp(p)}
\end{align*}
for $p=\sum_{\underline\alpha} p_{\underline\alpha}\underline x^{\underline\alpha}\in\mathbb{K}[\![\underline x]\!]$ and that
\[
\lead(P)=\{\lead(p)\,\vert\,p\in P\}
\]
for a subset $P\subset\mathbb{K}[\![\underline x]\!]$.
Note that the maximum exists by Dickson's Lemma \cite[Lem.~1.2.6]{GP02}.
The weighted degree $\deg_{\underline x}$ extends to $\mathbb{K}[\![\underline x]\!]$ by setting
\[
\deg_{\underline x}(p)=\deg\lexp(p)
\]
for $p\in\mathbb{K}[\![\underline x]\!]$.
Since $[\partial_i,x_i]=1$, $\deg_{\underline x}$ extends to $\mathbb{K}[\![\underline x]\!]\langle\underline\partial\rangle$ by setting
\[
\deg_{\underline x}(\partial_i)=-\deg_{\underline x}(x_i)>0.
\]

Let $\underline g=g_1,\dots,g_l$ be a standard basis of $\langle\underline f\rangle\mathbb{K}[\![\underline x]\!]$.
This means that $0\ne g_i\in\langle\underline f\rangle\mathbb{K}[\![\underline x]\!]$ and 
\begin{equation}\label{7}
\lead\bigl(\langle\underline f\rangle\mathbb{K}[\![\underline x]\!]\bigr)=\langle\lead(\underline g)\rangle\mathbb{K}[\![\underline x]\!]
\end{equation}
which implies that
\begin{equation}\label{8}
\langle\underline f\rangle\mathbb{K}[\![\underline x]\!]=\langle\underline g\rangle\mathbb{K}[\![\underline x]\!]
\end{equation}
by the division theorem.
Let
\[
\underline m=(m_i)_{i=1,\dots,\mu}=(\underline x^{\underline\beta})_{\underline x^{\underline\beta}\notin\langle\lead(\underline g)\rangle\mathbb{K}[\![\underline x]\!]}
\]
be increasingly ordered with respect to $<_{\underline x}$.
Then
\begin{equation}\label{6}
\mathbb{K}[\![\underline x]\!]=\langle\underline m\rangle\mathbb{K}\oplus\langle\lead(\underline g)\rangle\mathbb{K}[\![\underline x]\!]
\end{equation}
and hence, by (\ref{7}) and (\ref{8}), $\mathbb{K}[\![\underline x]\!]=\langle\underline m\rangle\mathbb{K}\oplus\langle\underline g\rangle\mathbb{K}[\![\underline x]\!]$.
Then $\underline m$ represents a $\mathbb{K}$--basis of 
\begin{equation}\label{5}
H/\langle\underline s\rangle H=\mathbb{K}[\![\underline x]\!]/\langle\underline g\rangle\mathbb{K}[\![\underline x]\!]=\langle\underline m\rangle\mathbb{K}
\end{equation}
and, by Nakayama's Lemma, $\underline m$ represents a minimal set of $\mathbb{K}[\![\underline s]\!]$--generators of $H$.
Note that if $H$ is free then it is free of rank $\mu$.
Let $U=(\overline u_i)_i\in\mathbb{K}[\![\underline x]\!]^{r\times l}$ be a matrix such that
\[
\underline g=\underline fU.
\]

\begin{rem}
If $\underline f\in\mathbb{K}[\underline x]$ then one can compute $\underline g$ and $U$ with coefficients in $\mathbb{K}[\underline x]$ by Lazard's method based on Buchberger's standard basis algorithm \cite[Lem.~1.7]{GP96} and homogenization.
In general, the power series $f_i\in\mathbb{K}[\![\underline x]\!]$ can be represented by generating functions $\SelectTips{cm}{}\xymatrix{\Nset^n\ar[r]&\mathbb{K}}$ and Buchberger's standard basis algorithm with respect to a local degree ordering computes such generating functions for $\underline g$ and $U$.
\end{rem}

Let $<_{\underline s}$ be a local degree ordering with respect to a weighted degree $\deg_{\underline s}$ on the monomials of $\mathbb{K}[\![\underline s]\!]$.
Let 
\[
<=(<_{\underline s},<_{\underline x})
\]
be the block ordering of $<_{\underline s}$ and $<_{\underline x}$ on the monomials of $\mathbb{K}[\![\underline s,\underline x]\!]$ and 
\[
\deg=\deg_{\underline s}+\deg_{\underline x}
\]
the sum of the weighted degrees $\deg_{\underline s}$ and $\deg_{\underline x}$.
This means that
\[
\underline s^{\underline\alpha'}\underline x^{\underline\beta'}<\underline s^{\underline\alpha''}\underline x^{\underline\beta''}\Leftrightarrow\underline s^{\underline\alpha'}<_{\underline s}\underline s^{\underline\alpha''}\vee\bigl(\underline s^{\underline\alpha'}=\underline s^{\underline\alpha''}\wedge\underline x^{\underline\beta'}<_{\underline x}\underline x^{\underline\beta''}\bigr)
\]
and
\[
\deg(\underline s^{\underline\alpha}\underline x^{\underline\beta})=\deg_{\underline s}(\underline s^{\underline\alpha})+\deg_{\underline x}(\underline x^{\underline\beta}).
\]
As before, we denote the leading exponent, resp. leading term, with respect to $<$ by $\lexp$, resp. $\lead$, and extend $\deg$ to $\mathbb{K}[\![\underline s,\underline x]\!]\langle\underline\partial\rangle$.
Note that $<$ is not a degree ordering with respect to $\deg$.
This means that
\[
\deg\lead\ne\deg.
\]
We denote by the leading exponent, resp. leading term, with respect to the partial ordering $<_{\underline s}$ on $\mathbb{K}[\![\underline s,\underline x]\!]$ by $\lexp_{\underline s}$, resp. $\lead_{\underline s}$, and the partial degree $<_{\underline s}$ on $\mathbb{K}[\![\underline s,\underline x]\!]$ by $\deg_{\underline s}$.
This means that
\begin{align*}
\lexp_{\underline s}(p)&=\max_{<_{\underline s}}\pi_{\Nset^m}(\supp(p)),\\
\lead_{\underline s}(p)&=p_{\lexp_{\underline s}(p)}\underline s^{\lexp_{\underline s}(p)},\\
\deg_{\underline s}(p)&=\deg\lexp_{\underline s}(p)
\end{align*}
for $p=\sum_\alpha p_{\underline\alpha}\underline s^{\underline\alpha}\in\mathbb{K}[\![\underline s,\underline x]\!]$ where $\SelectTips{cm}{}\xymatrix{\pi_{\Nset^m}:\Nset^m\times\Nset^n\ar[r] & \Nset^m}$ is the canonical projection.
We denote by $\min\deg$, resp. $\max\deg$, the minimal, resp. the maximal, degree of the components of a vector or a matrix.

Let
\[
\underline G=G_1,\dots,G_l=\underline FU=\underline g-\underline sDU.
\]
In the special fibre $\underline s=\underline0$, $\underline G$ induces the standard basis $\underline g$.
We call $\underline G$ a partial standard basis of the formal family $\underline F$.

The following example is taken from \cite[Sec.~8]{Sch02c}.

\begin{exmp}\label{27}
Let $\mathbb{K}=\Cset$, $\underline F=\underline\partial(f)-s\underline\partial$ as in Remark \ref{9}, and $f=x^5+x^2y^2+y^5$.
Note that $f$ defines a $T_{2,5,5}$ singularity at the origin.
Let $<_{(x,y)}$ be the local degree ordering with $\deg(x)=\deg(y)=-1$ and $x>y$.
Then one computes
\begin{align*}
\underline f&=2xy^2+5x^4,2x^2y+5y^4,\\
\underline g&=2x^2y+5y^4,2xy^2+5x^4,5x^5-5y^5,10y^6+25x^3y^4,\\
\mu&=11,\\
\underline m&=y^5,y^4,y^3,y^2,xy,y,x^4,x^3,x^2,x,1,\\
U&=\begin{pmatrix}
0 & 1 & x & -2xy\\
1 & 0 & -y & 2y^2+5x^3
\end{pmatrix}
\end{align*}
and hence
\begin{multline*}
\underline G=2x^2y+5y^4-s\partial_y,2xy^2+5x^4-s\partial_x,\\
5x^5-5y^5-sx\partial_x+sy\partial_y,10y^6+25x^3y^4+2sxy\partial_x-s(2y^2+5x^3)\partial_y.
\end{multline*}
\end{exmp}

We denote by $\underline F\underline x^{\Nset^n}=(F_i\underline x^{\underline\alpha})_{i,\underline\alpha}$ the generators of the $\Cset[\![\underline s]\!]$--module $\langle\underline F\rangle\mathbb{K}[\![\underline s,\underline x]\!]$.

\begin{lem}\label{10}
$H$ is a free $\mathbb{K}[\![s]\!]$--module if and only if $\underline G\underline x^{\Nset^n}$ is a standard basis of the $\mathbb{K}[\![\underline s]\!]$--module $\langle\underline F\rangle\mathbb{K}[\![\underline s,\underline x]\!]$.
\end{lem}
\begin{pf}
By (\ref{5}) and Nakayama's Lemma, $\underline m$ represents a minimal set of generators of $H$.
Since $H=\mathbb{K}[\![\underline s,\underline x]\!]/\langle\underline F\rangle\mathbb{K}[\![\underline s,\underline x]\!]$,
\[
\mathbb{K}[\![\underline s,\underline x]\!]=\langle\underline m\rangle\mathbb{K}[\![s]\!]+\langle\underline F\rangle\mathbb{K}[\![\underline s,\underline x]\!]
\]
and $H$ is free if and only if
\[
\langle\underline F\rangle\mathbb{K}[\![\underline s,\underline x]\!]\cap\langle\underline m\rangle\mathbb{K}[\![s]\!]=0.
\]
By (\ref{7}) and (\ref{6}), this is equivalent to
\[
\lead\bigl(\langle\underline F\rangle\mathbb{K}[\![\underline s,\underline x]\!]\bigr)=\langle\lead(\underline g)\rangle\mathbb{K}[\![\underline s,\underline x]\!]=\bigl\langle\lead\bigl(\underline G\underline x^{\Nset^n}\bigr)\bigr\rangle\mathbb{K}[\![\underline s]\!].\qed
\]
\end{pf}

By Proposition \ref{4}, 
\[
\Cset[\![s,\underline x]\!]/\langle\underline\partial(f)-s\underline\partial\rangle\Cset[\![s,\underline x]\!]\cong_{\Cset[\![s]\!]} H''
\]
is a free $\Cset[\![s]\!]$--module of rank $\mu$.
We shall now give an elementary proof of this fact.

\begin{prop}\label{21}
If $\mathbb{K}=\Cset$ and $\underline F=\underline\partial(f)-s\underline\partial$ as in Remark \ref{9} then $H$ is a free $\Cset[\![s]\!]$--module of rank $\mu$.
\end{prop}
\begin{pf}
Let $0\ne p\in\langle\underline m\rangle\Cset[\![s]\!]\cap\langle\underline\partial(f)-s\underline\partial\rangle\Cset[\![s,\underline x]\!]$.
Then $\lead(p)\in\langle\underline m\rangle\Cset[\![s]\!]$ and $p=(\underline\partial(f)+s\underline\partial)\overline q$ for some $\overline q\in\Cset[\![s,\underline x]\!]$ with maximal $\max\deg_s(\overline q)$.
By (\ref{7}) and (\ref{6}), this implies that $\underline\partial(f)\lead_s(\overline q)=0$ and hence, by Lemma \ref{2}, we may assume that there are $1\le i<j\le n$, $k\ge0$, and $r\in\Cset[\![\underline x]\!]$ such that
\[
\lead_s(\overline q)=s^kr\bigl(\partial_i(f)\overline e_j-\partial_j(f)\overline e_i\bigr).
\]
This implies that
\[
\underline\partial\lead_s(\overline q)=s^k\underline\partial(f)\bigl(\partial_j(r)\overline e_i-\partial_i(r)\overline e_j\bigr)
\]
and hence
\[
p=(\underline\partial(f)-s\underline\partial)\bigl(\overline q-\lead_s(\overline q)-s^{k+1}\bigl(\partial_j(r)\overline e_i-\partial_i(r)\overline e_j\bigr)\bigr).
\]
This is a contradiction to the maximality of $\max\deg_s(\overline q)$.
Hence,
\[
\langle\underline m\rangle\Cset[\![s]\!]\cap\langle\underline\partial(f)-s\underline\partial\rangle\Cset[\![s,\underline x]\!]=0
\]
and $H$ is free.\qed
\end{pf}

Our aim is now to define a filtration $V=(V_K)_{K\le0}$ on $\mathbb{K}[\![\underline s,\underline x]\!]$ by $\mathbb{K}[\![\underline s]\!]$--modules which is
\begin{enumerate}
\item a basis of the $\langle\underline s,\underline x\rangle$--adic topology on $\mathbb{K}[\![\underline s,\underline x]\!]$,
\item compatible with reduction with respect to the partial standard basis $\underline G$,
\item mapped by $\pi_H$ onto the basis $(\langle\underline s\rangle^KH)_{K\ge0}$ of the $\langle\underline s\rangle$--adic topology on $H$.
\end{enumerate}
This will lead to a normal form algorithm for $\underline H$.

For a given weighted degree $\deg_{\underline x}$, let the weighted degree $\deg_{\underline s}$ be such that
\begin{equation}\label{40}
\deg(s_j)\le\min\deg(\underline m)+\min\deg(\underline x)-\max\deg(\underline d^j).
\end{equation}
Let the strictly increasing sequence $N=(N_K)_{K\le0}$ be defined by 
\begin{equation}\label{41}
N_K=-K\min\deg(\underline s)-\min\deg(\underline x)+\max\deg(D).
\end{equation}
Let $V=(V_K)_{K\le0}$ be the strictly increasing filtration on $\mathbb{K}[\![\underline s,\underline x]\!]$ by $\mathbb{K}[\![\underline s]\!]$--modules 
\begin{equation}\label{42}
V_K=\bigl\{p\in\mathbb{K}[\![\underline s,\underline x]\!]\;\big\vert\deg(p)<N_K\bigr\}+\langle\underline s\rangle^{-K}\mathbb{K}[\![\underline s,\underline x]\!].
\end{equation}

\begin{rem}
For $\underline F=\underline\partial(f)-s\underline\partial$ as in Remark \ref{9}, we can choose
\begin{align*}
\deg(s)&=\min\deg(\underline m)+2\min\deg(\underline x),\\
N_K&=-K\deg(s)-2\min\deg(\underline x).
\end{align*}
\end{rem}

\begin{exmp}\label{28}
In example \ref{27}, $\deg(s)=-7$ and $N_K=7K+2$.
\end{exmp}

The following proposition is a generalization of \cite[Lem.~10]{Sch02c}.

\begin{prop}\label{11}\
\begin{enumerate}
\item $V=(V_K)_{K\le0}$ is a basis of the $\langle\underline s,\underline x\rangle$--adic topology.
\item If $\lead(\underline s^{\underline\alpha} G_k\underline x^{\underline\beta})\in V_K$ then also $\underline s^{\underline\alpha} G_k\underline x^{\underline\beta}\in V_K$.
\item $\pi_H(V_K)=\langle\underline s\rangle^{-K}H$.
\end{enumerate}
\end{prop}
\begin{pf}
\begin{enumerate}
\item This follows from (\ref{41}) and (\ref{42}).
\item Since $\underline g$ is a standard basis,
\begin{equation}\label{13}
\min\deg(\underline m)+\min\deg(\underline x)\le\min\deg(\underline g)
\end{equation}
and hence, by (\ref{40}),
\begin{align*}
\deg(\underline sD\overline u_k)&\le\max\{\deg(s_j\underline d^j\overline u_k)\mid 1\le j\le m\}\\
&\le\max\{\deg(s_j)+\max\deg(\underline d^j)+\max\deg(\overline u_k)\mid 1\le j\le m\}\\
&\le\max\{\deg(s_j)+\max\deg(\underline d^j)\mid 1\le j\le m\}\\
&\le\min\deg(\underline m)+\min\deg(\underline x)\\
&\le\min\deg(\underline g)\le\deg(g_k)\\
&=\deg\lead(g_k).
\end{align*}
Since $\underline s^{\underline\alpha} G_k\underline x^{\underline\beta}=\underline s^{\underline\alpha}(g_k-\underline sD\overline u_k)\underline x^{\underline\beta}$, this implies that
\[
\deg\lead(\underline s^{\underline\alpha} G_k\underline x^{\underline\beta})=\deg(\underline s^{\underline\alpha}\lead(g_k)\underline x^{\underline\beta})=\deg(\underline s^{\underline\alpha} G_k\underline x^{\underline\beta}).
\]
Hence, the claim follows from (\ref{42}).
\item Let $0\ne p\in V_K$ and $\underline s^{\underline\alpha} p_{\underline\alpha}=\lead_{\underline s}(p)$ with maximal $\vert\underline\alpha\vert<-K$ for fixed $p\mod\langle\underline F\rangle\mathbb{K}[\![\underline s,\underline x]\!]$.
Then, by (\ref{40}), 
\begin{align*}
\deg(p_{\underline\alpha})&=\deg(\underline s^{\underline\alpha} p_{\underline\alpha})-\deg(\underline s^{\underline\alpha})\\
&<-(K+\vert\underline\alpha\vert)\min\deg(\underline s)-\min\deg(\underline x)+\max\deg(D)\\
&\le\min\deg(\underline s)-\min\deg(\underline x)+\max\deg(D)\\
&\le\min\deg(\underline m)
\end{align*}
and hence, by (\ref{6}), $p_{\underline\alpha}\in\langle\underline g\rangle\mathbb{K}[\![\underline x]\!]$.
By the division theorem, there is a $\overline q\in\mathbb{K}[\![\underline x]\!]$ with $p_{\underline\alpha}=\underline g\overline q$ and $\lead(p_{\underline\alpha})\ge\lead(g_jq^j)$ for all $j$ and hence
\begin{equation}\label{15}
\max\deg(\overline q)\le\deg(p_{\underline\alpha})-\min\deg(\underline g).
\end{equation}
Then 
\begin{equation}\label{12}
p_{\underline\alpha}=\underline g\overline q=\underline fU\overline q\equiv\underline sDU\overline q\mod\langle\underline F\rangle\mathbb{K}[\![\underline s,\underline x]\!]
\end{equation}
and hence, by (\ref{40}), (\ref{13}), and (\ref{15}) 
\begin{align*}
\deg(\underline sDU\overline q)
&\le\max\deg(\underline sD)+\max\deg(U)+\max\deg(\overline q)\\
&\le\max\deg(\underline sD)+\max\deg(\overline q)\\
&\le\max\deg(\underline sD)-\min\deg(\underline g)+\deg(p_{\underline\alpha})\\
&\le\max\deg(\underline sD)-\min\deg(\underline m)-\min\deg(\underline x)+\deg(p_{\underline\alpha})\\
&\le\max\{\deg(s_j)+\max\deg(\underline d^j)\mid 1\le j\le m\}\\
&-\min\deg(\underline m)-\min\deg(\underline x)+\deg(p_{\underline\alpha})\\
&\le\deg(p_{\underline\alpha}).
\end{align*}
Hence, by (\ref{12}),
\[
p'=p-\lead_{\underline s}(p)+\underline s^{\underline\alpha}\underline sDU\overline q\equiv p\mod\langle\underline F\rangle\mathbb{K}[\![\underline s,\underline x]\!]
\]
with $\deg(p')\le\deg(p)<N_K$ and $\lead_{\underline s}(p')<_{\underline s}\lead_{\underline s}(p)$.
This contradicts to the maximality of $\vert\underline\alpha\vert$ and hence $p\in\langle\underline s\rangle^{-K}+\langle\underline F\rangle\mathbb{K}[\![\underline s,\underline x]\!]$.\qed
\end{enumerate}
\end{pf}

Proposition \ref{11} leads to the following normal form algorithm.

\begin{alg}\label{34}\
\begin{tabbing}
\quad\=\quad\=\quad\=\quad\kill
\textbf{proc} $\NF(p\in\mathbb{K}[\![\underline s,\underline x]\!],K\le0)$\\
\>\textbf{if} $p\in\langle\underline s\rangle^{-K}$ \textbf{then} $q:=p$\\
\>\textbf{else} \textbf{if} $\deg\lead(p)<N_K$ \textbf{or} $\lead_{\underline s}(p)\in\langle\underline s\rangle^{-K}$ \textbf{then} $q:=\lead_{\underline s}(p)$\\
\>\textbf{else} $q:=0$\\
\>$r:=p-q$\\
\>\textbf{if} $r=0$ \textbf{then} \textbf{return} $r\in\mathbb{K}[\![\underline s,\underline x]\!],\overline a\in\mathbb{K}[\![\underline s,\underline x]\!]^l,q\in\mathbb{K}[\![\underline s,\underline x]\!]$\\
\>\textbf{if} $\lead(r)\in\langle\lead(\underline g)\rangle$ \textbf{then}\\
\>\>$j:=\min\{i\mid\lead(r)\in\langle\lead(g_i)\rangle\}$\\
\>\>$r,\overline a,q':=\NF\bigl(r-\frac{\lead(r)}{\lead(g_j)}g_j-\underline sD\bigl(\frac{\lead(r)}{\lead(g_j)}\overline u_j\bigr),K\bigr)$\\
\>\>$\overline a:=\overline a+\frac{\lead(r)}{\lead(g_j)}\overline e_j$\\
\>\textbf{else}\\
\>\>$r',\overline a,q':=\NF(r-\lead(r),K)$\\
\>\>$r:=\lead(r)+r'$\\
\>$q:=q+q'$\\
\>\textbf{return} $r\in\mathbb{K}[\![\underline s,\underline x]\!],\overline a\in\mathbb{K}[\![\underline s,\underline x]\!]^l,q\in\mathbb{K}[\![\underline s,\underline x]\!]$.\\
\end{tabbing}
\end{alg}

The input of the algorithm $\NF$ is a power series $p\in\mathbb{K}[\![\underline s,\underline x]\!]$ and an integer $K\le0$, the output is a power series $r\in\mathbb{K}[\![\underline s,\underline x]\!]$, a column vector $\overline a$ with coefficients in $\mathbb{K}[\![\underline s,\underline x]\!]$, and a power series $q\in\mathbb{K}[\![\underline s,\underline x]\!]$.
We denote the components of $\NF$ by
\[
(\NF_1(p,K),\NF_2(p,K),\NF_3(p,K))=(r,\overline a,q)=\NF(p,K)
\]
for $p\in\mathbb{K}[\![\underline s,\underline x]\!]$ and $K\le0$.

\begin{exmp}\label{29}
In example \ref{27} using \ref{28}, one computes
\[
\NF_1(f\underline m,-2)=\underline m(A_0+sA_1)
\]
where $A_0,A_1\in\Cset^{11\times11}$ such that
\[
\setcounter{MaxMatrixCols}{20}
A_0+sA_1=
\begin{pmatrix}
\frac{3}{2}s&0&0&0&-\frac{25}{4}s&0&0&0&0&0&-\frac{1}{2}\\
0&\frac{13}{10}s&0&0&0&0&0&0&0&-\frac{75}{16}s&0\\
0&0&\frac{11}{10}s&0&0&0&0&0&-\frac{1}{4}s&0&0\\
0&0&0&\frac{9}{10}s&0&0&0&0&0&0&0\\
0&0&0&0&s&0&0&0&0&0&0\\
0&0&0&0&0&\frac{7}{10}s&0&0&0&0&0\\
0&0&0&0&0&-\frac{75}{16}s&\frac{13}{10}s&0&0&0&0\\
0&0&0&-\frac{1}{4}s&0&0&0&\frac{11}{10}s&0&0&0\\
0&0&0&0&0&0&0&0&\frac{9}{10}s&0&0\\
0&0&0&0&0&0&0&0&0&\frac{7}{10}s&0\\
0&0&0&0&0&0&0&0&0&0&\frac{1}{2}s
\end{pmatrix}.
\]
\end{exmp}

Figure \ref{17} illustrates a reduction step in $\NF$.
The $1$--dimensional $\mathbb{K}$--vector space spanned by a monomial in $\langle\underline m\rangle\mathbb{K}[\![\underline s]\!]$, resp. in $\langle\lead(\underline g)\rangle\mathbb{K}[\![\underline s,\underline x]\!]$, is depicted by a big, resp. small, bullet.
The monomial at the tail of the arrow is replaced by a power series with support above the dotted line meeting the head of the arrow.
The $\mathbb{K}[\![s]\!]$--submodule $V_K$ generated by the monomials above the dashed line is invariant with respect to such a reduction step.

\begin{figure}[ht]
\caption{A reduction step in $\NF$}\label{17}
\setlength{\unitlength}{1mm}
\begin{picture}(130,70)(-5,0)
\path(0,70)(0,10)(130,10)
\put(0,70){\path(0,0)(-1,-2)\path(0,0)(1,-2)}
\put(130,10){\path(0,0)(-2,-1)\path(0,0)(-2,1)}
\put(-2,70){\makebox(0,0)[rt]{$\underline{s}$}}
\put(130,8){\makebox(0,0)[rt]{$\underline{x}$}}
\matrixput(0,10)(5,0){5}(0,5){12}{\circle*{1}}
\matrixput(25,10)(5,0){21}(0,5){12}{\circle*{.5}}
\put(0,60){\path(-1,0)(1,0)\put(-1,0){\makebox(0,0)[r]{$K$}}}
\dashline[1]{1}(0,60)(20,60)(120,10)
\put(72,35){\makebox(0,0)[bl]{$V_{-K}$}}
\dashline[1]{1}(23,70)(23,0)(130,0)
\put(24,1){\makebox(0,0)[bl]{$\langle\lead(\underline g)\rangle\mathbb{K}[\![\underline s,\underline x]\!]$}}
\dashline[1]{1}(22,70)(22,0)(0,0)
\put(21,1){\makebox(0,0)[br]{$\langle\underline m\rangle\mathbb{K}[\![\underline s]\!]$}}
\dottedline{1}(0,55)(60,25)(65,25)(75,20)(100,20)(120,10)
\path(60,25)(70,20)
\put(60,25){\path(0,0)(2,-0.5)\path(0,0)(1.5,-1.5)}
\end{picture}
\end{figure}

\begin{lem}
$\NF$ terminates.
\end{lem}
\begin{pf}
For fixed leading exponent $\lexp_{\underline s}(p)$ with respect to $\underline s$, the leading term $\lead(p)$ is strictly decreasing with weighted degree $\deg\lead(p)\ge N_K$.
Since there are only finitely many monomials with fixed $\lexp_{\underline s}$ and $\deg\ge N_K$, $\lexp_{\underline s}(p)$ decreases after finitely many steps.
Since $<=(<_{\underline s},<_{\underline x})$ is a block ordering and $<_{\underline s}$ is a degree ordering, this implies that $p\in\langle\underline s\rangle^{-K}$ after finitely many steps.\qed
\end{pf}

The following lemma states that $\NF_1(\cdot,L)$ is a reduced normal form modulo $V_L$ with $\NF_1(V_K,L)\subset V_K$ for $L<K\le0$.

\begin{lem}\label{18}
Let $L<K\le0$, $p\in V_K$, and $(r,\overline a,q)=\NF(p,L)$.
Then 
\begin{enumerate}
\item $p=\underline G\overline a+r+q$,
\item $a^i\in\mathbb{K}[\underline s,\underline x]$ with $\lead(G_ia^i)\le\lead(p)$ for $i=1,\dots,l$,
\item $r\in\langle\underline m\rangle\bigoplus_{L<\vert\underline\alpha\vert\le K}\mathbb{K}\underline s^{\underline\alpha}$ with $\lead(r)\le\lead(p)$,
\item $q\in V_L$.
\item If $p\in\langle\underline m\rangle\bigoplus_{L<\vert\underline\alpha\vert\le K}\mathbb{K}\underline s^{\underline\alpha}$ then $(r,\overline a,q)=(p,\overline0,0)$.
\item If $p\equiv p'\mod V_L$ then $r=r'$. 
\end{enumerate}
\end{lem}
\begin{pf}
By Proposition \ref{11}.2, $\NF$ preserves the condition $p\in V_K$.
Hence, the claim follows immediately from the definition of $\NF$.\qed
\end{pf}

By Proposition \ref{11}.1, V is a basis of the $\langle\underline s,\underline x\rangle$--adic topology and, by Lemma \ref{18}, $\NF_1(V_K,L)\subset V_K$ and $\NF_3(V_K,L)\subset V_L$ for $L<K\le0$.
Since $\mathbb{K}[\![\underline s,\underline x]\!]$ is complete with respect to the $\langle\underline s,\underline x\rangle$--adic topology, 
\[
\mathbb{K}[\![\underline s,\underline x]\!]=\underset{\underset{K}{\longleftarrow}}{\lim}(\mathbb{K}[\![\underline s,\underline x]\!]/V_K)
\]
and hence $\NF$ induces a reduced normal form on $\mathbb{K}[\![\underline s,\underline x]\!]$ as follows.

\begin{defn}
Let $K=(K_i)_{i\ge 0}$ be a strictly decreasing sequence and 
\[
\NF(p)=(\NF_1(p),\NF_2(p))=\Bigl(\sum_{i\ge0}r_i,\sum_{i\ge0}\overline a_i\Bigr)
\]
for $p\in\mathbb{K}[\![\underline s,\underline x]\!]$ where $p_0=p$ and $r_i,\overline a_i,p_{i+1}=\NF(p_i,K_i)$ for $i\ge0$.
\end{defn}

Note that $\NF$ depends on the choice of the sequence $K$.

\begin{lem}\label{19}
Let $(r,\overline a)=\NF(p)$.
Then
\begin{enumerate}
\item $p=\underline G\overline a+r$,
\item $a^i\in\mathbb{K}[\![\underline s,\underline x]\!]$ with $\lead(G_ia^i)\le\lead(p)$ for $i=1,\dots,l$,
\item $r\in\langle\underline m\rangle\mathbb{K}[\![\underline s]\!]$ with $\lead(r)\le\lead(p)$.
\item If $p\in\langle\underline m\rangle\mathbb{K}[\![\underline s]\!]$ then $(r,\overline a)=(p,\overline0)$.
\end{enumerate}
\end{lem}
\begin{pf}
This follows immediately from Proposition \ref{11}.1 and Lemma \ref{18}.\qed
\end{pf}

The following proposition describes $\NF_1$ as a map of $\Cset[\![\underline s]\!]$--modules.

\begin{prop}\label{20}\
$\NF_1$ is a $\mathbb{K}[\![\underline s]\!]$--linear map 
\[
\SelectTips{cm}{}\xymatrix@C=20pt{
\mathbb{K}[\![\underline s,\underline x]\!]\ar@{->>}[rd]_-{\pi_H}\ar@{->>}[rr]^-{\NF_1} && \langle\underline m\rangle\mathbb{K}[\![\underline s]\!]\ar@{->>}[ld]^-{\pi_H}\\
& H}
\]
with $\NF_1^2=\NF_1$.
In particular, if $H$ is a free $\mathbb{K}[\![\underline s]\!]$--module then $\NF_1$ induces the $\mathbb{K}[\![\underline s]\!]$--section
\[
\SelectTips{cm}{}\xymatrix@C=40pt{\mathbb{K}[\![\underline s,\underline x]\!]\ar@{->>}[r]_-{\pi_H} & H\ar@/_1pc/[l]_-{\NF_1}}
\]
of the canonical projection $\pi_H$ with image $\NF_1(H)=\langle\underline m\rangle\mathbb{K}[\![\underline s]\!]$.
This means that $\NF_1$ is the $\underline m$--basis representation.
\end{prop}
\begin{pf}
By definition of $\NF$, $\NF_1$ is $\mathbb{K}[\![\underline s]\!]$--linear.
By Lemma \ref{19}, $\NF_1$ is a map over $H$.
By Lemma \ref{19}.4, $\NF_1$ is the identity on its image $\langle\underline m\rangle\mathbb{K}[\![\underline s]\!]$.\qed
\end{pf}

If $H$ is not a free $\mathbb{K}[\![\underline s]\!]$--module then its relations can be computed as follows.
By Proposition \ref{20},
\begin{align*}
H&\cong_{\mathbb{K}[\![\underline s]\!]}\langle\underline m\rangle\mathbb{K}[\![\underline s]\!]\big/\bigl(\langle\underline m\rangle\mathbb{K}[\![\underline s]\!]\cap\langle\underline F\rangle\mathbb{K}[\![\underline s,\underline x]\!]\bigr)\\
&=\langle\underline m\rangle\mathbb{K}[\![\underline s]\!]\big/\NF_1\bigl(\langle\underline F\rangle\mathbb{K}[\![\underline s,\underline x]\!]\bigr)\\
&=\langle\underline m\rangle\mathbb{K}[\![\underline s]\!]\big/\bigl\langle\NF_1\bigl(\underline F\underline x^{\Nset^n}\bigr)\bigr\rangle\mathbb{K}[\![\underline s]\!]
\end{align*}
and, in particular, by Lemma \ref{18}.3,
\[
H/\langle\underline s\rangle^{-K}H\cong_{\mathbb{K}[\![\underline s]\!]}\langle\underline m\rangle\mathbb{K}[\![\underline s]\!]\big/\bigl(\bigl\langle\NF_1\bigl(\underline F\underline x^{\Nset^n}\backslash V_K,K\bigr)\bigr\rangle\mathbb{K}[\![\underline s]\!]+\langle\underline m\rangle\langle\underline s\rangle^{-K}\mathbb{K}[\![\underline s]\!]\bigr)
\]
where $\underline F\underline x^{\Nset^n}\backslash V_K$ is a finite subset of $\mathbb{K}[\![\underline s,\underline x]\!]$.

Finally, we return to our starting point.
Let $\mathbb{K}=\Cset$ and $\underline F=\underline\partial(f)-s\underline\partial$ as in Remark \ref{9}.
Then, by Proposition \ref{4}, $H\cong_{\Cset[\![s]\!]}\widehat H''$ is the formal Brieskorn lattice and, by Proposition \ref{21}, $H$ is a free $\Cset[\![s]\!]$--module of rank $\mu$.
We define the matrix $A^{\underline m}\in\Cset[\![s]\!]^{\mu\times\mu}$ by 
\[
\underline mA^{\underline m}=t\underline m.
\]
Then, by (\ref{24}), $A^{\underline m}+s^2\partial_s$ is the $\underline m$--basis representation of $t$.
This means that there is a commutative diagram
\[
\SelectTips{cm}{}\xymatrix@C=40pt{
\widehat H''\ar[r]^-t & \widehat H''\\
\Cset[\![\underline s]\!]^\mu\ar[u]_-\sim^-{\underline m}\ar[r]^-{A^{\underline m}+s^2\partial_s} & \Cset[\![\underline s]\!]^\mu\ar[u]_-\sim^-{\underline m}.}
\]
By (\ref{22}) and propsosition \ref{20},
\[
\underline mA^{\underline m}=\NF_1(f\underline m).
\]

\begin{exmp}
In example \ref{27} using \ref{29},
\[
A^{\underline m}\equiv A_0+sA_1\mod s^2\Cset[\![s]\!]^{11\times11}.
\]
By \cite[Sec.~8]{Sch02c}, the non--diagonal terms of $A_1$ and the terms of $A^{\underline m}$ in $s^2\Cset[\![s]\!]^{11\times11}$ can be eliminated by transforming $\underline m$ to a good $\Cset[\![s]\!]$--basis of $\widehat H''$.
Then
\[
A^{\underline m}=A_0+sA_1
\]
and $A_0$ and $A_1$ represent M.~Saito's endomorphisms \cite{Sai89}.
\end{exmp}

All the algorithms in \cite{Sch01,SS01,Sch02a,Sch02c} require the computation of the matrix $A^{\underline m}$ for a $\Cset[\![s]\!]$--basis $\underline m$ of $\widehat H''$.

\section{Examples and timings}

Algorithm \ref{34} for the case of Remark \ref{9} and the algorithms in \cite{Sch02a,Sch02c,Sch03d} are implemented in the {\sc Singular} \cite{GPS03} library {\tt gmssing.lib} \cite{Sch03c}.
We use this implementation on a {\sc Pentium\,III\,M 1\,GHz} machine with $512$\,MB of memory plus $1$\,GB of swap memory.
For several polynomials $f\in\Cset[\underline x]$ with isolated critical point at the origin, we compute
\begin{enumerate}
\item the local Bernstein-Sato polynomial,
\item the spectral pairs, and
\item M.~Saito's endomorphisms $A_0$ and $A_1$.
\end{enumerate}
For $i=1,2,3$, we denote by $t_i$ the corresponding computation time in seconds and by $K_i$ the maximal $K$ occurring in $\NF$ during the computation.
By $t_K$ we denote the time in seconds needed to compute $A^{\underline m}\mod s^K\Cset[\![s]\!]^{\mu\times\mu}$.
All computation times are rounded off.

The local Bernstein-Sato polynomial at the origin for the examples in \cite[Tab.~1]{Nor02} can be computed, each in less than one second.
Table \ref{33} shows the timings for the examples in \cite[Tab.~2]{Nor02}.
By \cite{Sch02a}, it suffices to compute $A^{\underline m}\mod s^{K_0}\Cset[\![s]\!]^{\mu\times\mu}$ where $K_0=2(\mu+n-1)$ in order to compute all of the above invariants.
Table \ref{30} shows the results for the examples in \cite[Tab.~2]{Sch01} and that this a priori bound is useless in practice.

\begin{table}[ht]
\caption{Local Bernstein-Sato polynomial $b$ for $f=x^{n_1}+y^{n_2}+z^{n_3}+x^{m_1}y^{m_2}z^{m_3}$}\label{33}
\begin{tabular}{|c||c|c|c|c|c|c|c|c|}
\hline
$\underline n$ & $6,6,6$ & $7,7,7$ & $7,7,7$ & $9,9,9$ & $6,6,7$ & $6,6,7$ & $6,6,7$ & $6,7,7$\\
\hline
$\underline m$ & $4,4,4$ & $2,2,2$ & $2,2,3$ & $3,3,3$ & $2,2,2$ & $3,3,3$ & $4,4,4$ & $2,2,2$\\
\hline
$\mu$ & $125$ & $167$ & $216$ & $512$ & $138$ & $150$ & $150$ & $152$\\
\hline
$\deg(b)$ & $13$ & $28$ & $17$ & $23$ & $60$ & $48$ & $52$ & $62$\\
\hline
$K_1$ & $3$ & $4$ & $2$ & $2$ & $4$ & $3$ & $3$ & $4$\\
\hline
$t_1$ & $1$ & $5$ & $2$ & $50$ & $6$ & $1$ & $2$ & $6$\\
\hline
\end{tabular}
\end{table}

\begin{table}[ht]
\caption{Spectral pairs and M.~Saito's endomorphisms for \cite[Tab.~2]{Sch01}}\label{30}
\begin{tabular}{|c|c|c|c|c|c|c|c|c|c|c|}
\hline
 & $f$ & $\mu$ & $K_2$ & $t_{K_2}$ & $t_2$ & $K_3$ & $t_{K_3}$ & $t_3$ & $K_0$ & $t_{K_0}$\\
\hline
\hline
$Z_{1,1}$ & 
$x^{3}y+x^{2}y^{3}+y^{8}$ & 
$16$ &
$2$ & $0$ & $0$ & $4$ & $0$ & $10$ & $34$ & $180$\\
\hline
$W_{1,1}$ & 
$x^{4}+x^{2}y^{3}+y^{7}$ & 
$16$ &
$2$ & $0$ & $0$ & $4$ & $0$ & $2$ & $34$ & $39$\\
\hline
$W^\#_{1,1}$ & 
$x^{4}+2x^{2}y^{3}+xy^{5}+y^{6}$ & 
$16$ &
$2$ & $0$ & $1$ & $4$ & $1$ & $54$ & $34$ & $495$\\
\hline
$Q_{2,1}$ & 
$x^{3}+yz^{2}+x^{2}y^{2}+y^{7}$ & 
$15$ &
$2$ & $0$ & $1$ & $4$ & $0$ & $4$ & $34$ & $58$\\
\hline
$Q_{2,2}$ & 
$x^{3}+yz^{2}+x^{2}y^{2}+y^{8}$ & 
$16$ &
$2$ & $0$ & $0$ & $4$ & $0$ & $1$ & $36$ & $20$\\
\hline
$S_{1,1}$ & 
$x^{2}z+yz^{2}+x^{2}y^{2}+y^{6}$ & 
$15$ &
$2$ & $0$ & $0$ & $4$ & $0$ & $1$ & $34$ & $43$\\
\hline
$S_{1,2}$ & 
$x^{2}z+yz^{2}+x^{2}y^{2}+y^{7}$ & 
$16$ &
$2$ & $0$ & $0$ & $4$ & $0$ & $0$ & $36$ & $27$\\
\hline
$S^\#_{1,1}$ & 
$x^{2}z+yz^{2}+y^{3}z+xy^{4}$ & 
$15$ &
$2$ & $0$ & $1$ & $4$ & $0$ & $18$ & $34$ & $687$\\
\hline
$S^\#_{1,2}$ & 
$x^{2}z+yz^{2}+y^{3}z+x^{2}y^{3}$ & 
$16$ &
$2$ & $0$ & $1$ & $4$ & $0$ & $2$ & $36$ & $252$\\
\hline
$U_{1,1}$ & 
$x^{3}+xz^{2}+xy^{3}+y^{2}z^{2}$ & 
$15$ &
$2$ & $0$ & $1$ & $4$ & $0$ & $2$ & $34$ & $66$\\
\hline
$U_{1,2}$ & 
$x^{3}+xz^{2}+xy^{3}+y^{4}z$ & 
$16$ &
$2$ & $0$ & $0$ & $4$ & $0$ & $0$ & $36$ & $12$\\
\hline
$V_{1,1}$ & 
$x^{2}y+z^{2}y^{2}+z^{4}+y^{5}$ & 
$16$ &
$2$ & $0$ & $0$ & $6$ & $0$ & $2$ & $36$ & $12$\\
\hline
$V^\#_{1,1}$ & 
$x^{2}y+y^{4}+xz^{3}+yz^{3}$ & 
$16$ &
$2$ & $0$ & $1$ & $6$ & $0$ & $1$ & $36$ & $10$\\
\hline
\end{tabular}
\end{table}

By Remark \ref{9}, (\ref{8}) and (\ref{5}), the coefficient of $s^K$ in the $\underline m$--basis representation in $\widehat H''$ is defined by a division by the ideal $\langle\underline\partial(f)\rangle\Cset[\![\underline x]\!]$ where the output of the division for $s^K$ defines the input for the division for $s^{K+1}$.
Therefore, the complexity of the data and the computation time increases rapidly with $K$.
In \cite[Sec.~10.2]{Sch01}, we compute such a division by a sequence of weak normal form computations.
Table \ref{32} shows the time needed to compute $A^{\underline m}\mod s^K\Cset[\![s]\!]^{\mu\times\mu}$ for example \ref{27} and increasing $K$ using this method and the {\sc Singular} command {\tt division}.
The computation fails in degree $K=8$ after more than one hour due to lack of memory.
In the algorithm $\NF$, the above sequence of full divisions is replaced by a sequence of partial divisions.
Table \ref{31} shows the time needed to compute the same result using the algorithm $\NF$.
The situation is similar for other examples.

\begin{table}[ht]
\caption{$A^{\underline m}$ for $f=x^5+x^2y^2+y^5$ using {\tt division}}\label{32}
\begin{tabular}{|c||c|c|c|c|c|}
\hline
$K$ & $4$ & $5$ & $6$ & $7$ & $8$\\
\hline
$t_K$ & $0$ & $1$ & $10$ & $283$ & $\infty$\\
\hline
\end{tabular}
\end{table}

\begin{table}[ht]
\caption{$A^{\underline m}$ for $f=x^5+x^2y^2+y^5$ using $\NF$}\label{31}
\begin{tabular}{|c||c|c|c|c|c|c|c|c|c|c|}
\hline
$K$ & $20$ & $40$ & $60$ & $80$ & $100$ & $120$ & $140$ & $160$ & $180$ & $200$\\
\hline
$t_K$ & $1$ & $2$ & $8$ & $18$ & $35$ & $64$ & $106$ & $165$ & $247$ & $354$\\
\hline
\end{tabular}
\end{table}

\nocite{Sch03a}
\bibliographystyle{elsart-num}
\bibliography{nfbl}

\end{document}